\documentclass[a4paper, 10pt, twoside, notitlepage]{amsart}

\usepackage{amsmath,amscd}
\usepackage{amssymb}
\usepackage{amsthm}
\usepackage{comment}
\usepackage{graphicx}
\usepackage{epstopdf} 
\usepackage{mathrsfs}
\usepackage[ocgcolorlinks, linkcolor=blue]{hyperref}

\usepackage[ocgcolorlinks,linkcolor=blue]{hyperref}

\newcommand{\LC}{\left(}
\newcommand{\RC}{\right)}

\theoremstyle{plain}
\newtheorem{thm}{Theorem}[section]
\newtheorem{prop}{Proposition}[section]

\newtheorem{lem}[prop]{Lemma}

\newtheorem{assum}[prop]{Assumption}

\newtheorem{rmk}[prop]{Remark}

\numberwithin{equation}{section}
\newcommand {\R} {\mathbb{R}} 
 \newcommand {\N} {\mathbb{N}}
 
\newcommand {\p} {\partial}

\newcommand {\veps} {\varepsilon}

\newcommand {\supp} {\text{supp}}
\newcommand{\tre}{\textcolor{black}}

\newcommand{\norm}[1]{\lVert #1 \rVert}

\newcommand{\dx}[1][x]{\ensuremath{\,{\rm{d}} #1}}
\newcommand\mycom[2]{\genfrac{}{}{0pt}{}{#1}{#2}} 
\newcommand{\kommentar}[1]{}

\title[Simultaneous recovery semilinear elliptic inverse problems]{Simultaneous recovery of piecewise analytic coefficients in a semilinear elliptic equation}

\author[B. Harrach]{Bastian Harrach}
\address{Institute for Mathematics, Goethe University Frankfurt, Germany}
\email{harrach@math.uni-frankfurt.de}

\author[Y.-H. Lin]{Yi-Hsuan Lin}
\address{Department of Applied Mathematics, National Yang Ming Chiao Tung University, Taiwan}
\email{yihsuanlin3@gmail.com}

\begin{document}
	
	\maketitle

	\begin{abstract}
	In this short note, we investigate simultaneous recovery inverse problems for semilinear elliptic equations with partial data. The main technique is based on higher order linearization and monotonicity approaches. With these methods at hand, we can determine the diffusion and absorption coefficients together with the shape of a cavity simultaneously by knowing the corresponding localized Dirichlet-Neumann operator.

		\medskip
		
		
		
	\end{abstract}

	
	\let\thefootnote\relax\footnotetext{This is the accepted version of an article published in \emph{Nonlinear Analysis} \textbf{228}, 113188, 2023 (\url{https://doi.org/10.1016/j.na.2022.113188}). 

This accepted manuscript version is made available under the CC-BY-NC-ND 4.0 license
(\url{https://creativecommons.org/licenses/by-nc-nd/4.0/}).}

	
	

	\section{Introduction}
	In this note, we investigate some inverse problems for semilinear elliptic equations. 
	Inverse problems for nonlinear partial differential equations (PDEs) have been paid a lot of \tre{attention} in the past few decades. The main method to study inverse problems for nonlinear PDEs \tre{relies} on suitable linearization techniques, and the linearization approaches can be traced back to the pioneer work by Isakov \cite{isakov1993uniqueness}. In \cite{isakov1993uniqueness}, he demonstrated that the first linearization of the corresponding \tre{Dirichlet-to-Neumann} (DN) map of the semilinear parabolic equation agrees to the DN map of the associated linearized equation. Hence, related known results \tre{on} inverse boundary value problems for linear equations can be expected to apply, such that one is able to solve inverse problems for the nonlinear equations.

	For the semilinear elliptic equation $\Delta u +a(x,u)=0$ in a domain, the inverse problem of determining the coefficient $a(x,u)$ was treated in \cite{isakov1994global, sun2010inverse}, for $n\geq 3$, and in \cite{victorN, sun2010inverse, imanuvilovyamamoto_semilinear} for $n=2$.  In addition, for quasilinear elliptic equations, related inverse boundary value problems have also been studied by  \cite{sun1996, sun1997inverse, kang2002identification, liwang_navierstokes, munoz2020calderon,kian2020partial,carstea2021calder}. Meanwhile, researchers also worked on inverse problems for the degenerate $p$-Laplace equation~\cite{salo2012inverse, brander2018monotonicity}, and for the fractional semilinear Schr\"odinger equation \cite{lai2019global,LL2020inverse,lin2020monotonicity}. The interior unique determination problem for quasilinear equations on Riemannian manifolds was recently studied in~\cite[Section 6]{lassas2018poisson} \tre{via the source-to-solution map}. We \tre{also} refer the readers to \cite{sun2005, uhlmann2009calderon} for more introduction and discussions on \tre{related} inverse problems for nonlinear elliptic equations.

	In order to show the results in this work, we utilize a method that has been \tre{introduced by \cite{KLU2018inverse} for nonlinear hyperbolic equations} and developed in \cite{FO2019semilinear,LLLS2019nonlinear,LLLS2019partial,LLST2020nonlinear_fractional}. The method is called the \emph{higher order linearization}, which introduces particular parameters to reduce a semilinear elliptic equation into different linearized elliptic equations. In \cite{FO2019semilinear,LLLS2019nonlinear}, the authors studied inverse boundary value problem with full boundary measurements. In addition, simultaneous recovery inverse problems for semilinear PDEs were also considered in \cite{LLLS2019partial,LLST2020nonlinear_fractional,lin2021determining,lin2021simultaneous}, and \cite{KU2019partial,LLLS2019partial,KU2019derivative_partial} studied the Calder\'on type problem with partial data independently. A similar approach was utilized to study some related inverse problems for fractional semilinear elliptic equations \cite{LL2020inverse,lin2020monotonicity}. It can be noted that uniqueness results for the coefficients of the nonlinear term are easier to obtain than the uniqueness result for the coefficients of the linear term.
		
	In this note, we will not use complex geometrical optics solutions as in several of the previously cited works. Instead, after using the afore-mentioned \emph{higher order linearization} approach, we treat the resulting linearized elliptic equations by a combination of the \emph{monotonicity method} and \emph{localized potentials}. This line of reasoning was initiated by \cite{gebauer2008localized} and applied to various inverse problems, such as  \cite{harrach2009uniqueness,harrach2010exact,harrach2012simultaneous,arnold2013unique,harrach2013monotonicity,barth2017detecting,harrach2017local,brander2018monotonicity,griesmaier2018monotonicity,harrach2018helmholtz,harrach2018localizing,seo2018learning,harrach2019nonlocal-mono1,harrach2019dimension}. There are also related works on practical reconstruction methods based on monotonicity properties \cite{tamburrino2002new,harrach2015combining,harrach2015resolution,harrach2016enhancing,maffucci2016novel,tamburrino2016monotonicity,garde2017comparison,garde2017convergence,su2017monotonicity,ventre2017design,harrach2018monotonicity,zhou2018monotonicity,garde2019regularized,eberle2021shape,eberle2022monotonicity}.

 We next introduce the mathematical model in this work.	Let $\Omega \subset \R^n$ be a bounded connected domain with a $C^\infty$-smooth boundary $\p \Omega$, for $n\geq 2$, and $D\Subset \Omega$ is an open subset \tre{with a $C^\infty$-smooth boundary $\p D$} such that $\Omega \setminus \overline{D}$ is connected. Let us consider the following semilinear elliptic equation with diffusion 
\begin{align}\label{Main equation_obstacle}
	\begin{cases}
		-\nabla \cdot (\sigma(x) \nabla u) +a(x,u)=0 &\text{ in }\Omega\setminus \overline{D}, \\
		u=0& \text{ on }\p D, \\
		u=f & \text{ on }\p \Omega.
	\end{cases}
\end{align}
\tre{In particular, when $D=\emptyset$, the equation \eqref{Main equation_obstacle} becomes 
	\begin{align}\label{Main equation}
	\begin{cases}
		-\nabla \cdot (\sigma(x) \nabla u) +a(x,u)=0 &\text{ in }\Omega, \\
		u=f & \text{ on }\p \Omega.
	\end{cases}
\end{align}
}
	In this work, we first prove a local well-posedness result for  \eqref{Main equation_obstacle} under the assumption that the diffusivity $\sigma=\sigma(x)$, and the power series terms of the lower order coefficient $a=a(x,y)$, are \emph{piecewise real-analytic}, see Assumption~\ref{assum:pcw_anal} for the precise definition. 
Throughout this paper, we also synonymously use the term "analytic" instead of "real-analytic" for the sake of brevity. Note that a local well-posedness of \eqref{Main equation_obstacle} was already shown in \cite{FO2019semilinear,KU2019derivative_partial,LLLS2019nonlinear} under different regularity assumptions. 
Our case differs from previous results as we allow discontinuities of the coefficients, and we therefore give a detailed proof for the well-posedness result.

 To summarize the (local) well-posedness result, let $\Gamma \subset\p \Omega$ be a relatively open subset, and $0<\alpha<1$. We prove that there exists $\varepsilon>0$ and $C>0$ such that for all Dirichlet boundary data
\[
f\in N_\varepsilon:=\left\{f\in C^{\alpha}_0(\Gamma): \, \|f \|_{C^{\alpha}(\Gamma)} \leq \varepsilon \right\},
\]
there is a unique solution $u=u_f \in  V$ of \eqref{Main equation} that also satisfies
\[
\norm{u}_{V}\leq C \tre{ \varepsilon},
\] 
where 
\begin{align}\label{Function space V}
	V:=\left\{ v\in H^1(\Omega):\ -\nabla\cdot (\sigma\nabla v)\in L^\infty(\Omega),\ v|_{\partial \Omega}\in C^{\alpha}_0(\Gamma)\right\},
\end{align}
and 
\begin{align}\label{norm V}
	\norm{v}_V:=\norm{v}_{H^1(\Omega)}+\norm{\nabla\cdot (\sigma\nabla v)}_{L^\infty(\Omega)}+\norm{v|_{\partial \Omega}}_{C^{\alpha}(\Gamma)}.
\end{align}
Moreover, the local well-posedness also holds for \eqref{Main equation_obstacle}, when the domain $\Omega$ in \eqref{Function space V} and \eqref{norm V} is replaced by $\Omega \setminus \overline{D}$.

With these results at hand, one can define the corresponding (partial) DN operator
	\begin{align*}
	\Lambda_{\sigma,a}^{\Gamma}: N_\varepsilon &\to H^{-1/2}(\Gamma), \qquad \Lambda_{\sigma,a}^{\Gamma}(f) :=\left. \sigma\p_{\nu}u_f \right| _{\Gamma},
	\end{align*}
     for some sufficiently small number $\varepsilon>0$, where $u_f$ is the unique solution of \eqref{Main equation} and $\nu$ is the unit outer normal on $\p \Omega$. Likewise, one can also define the DN operator 
\begin{align*}
	\Lambda_{\sigma,a,D}^{\Gamma}: N_\varepsilon &\to H^{-1/2}(\Gamma), \qquad \Lambda_{\sigma,a,D}^{\Gamma}(f) :=\left. \sigma\p_{\nu}u_f \right| _{\Gamma},
\end{align*}	
  for some sufficiently small number $\varepsilon>0$, where $u_f$ is the unique solution of \eqref{Main equation_obstacle} and $\nu$ is the unit outer normal on $\p \Omega$. We then study the following simultaneous recovery inverse problems:
  \begin{itemize}
  	\item[(1)]  Can one simultaneously identify $\sigma$ and $a$ by knowing the partial measurements $\Lambda_{\sigma,a}^{\Gamma}$?
  	
  	\item[(2)]  Can one simultaneously identify $\sigma$, $a$ and $D$ by knowing the partial measurements $\Lambda_{\sigma,a,D}^{\Gamma}$?
   \end{itemize}
	We will give affirmative answers to both questions in this paper.

	The paper is structured as follows. In Section \ref{Sec:2}, we state our main results of this note, and the proofs are given in Section \ref{Sec:3}. The main methods depend on suitable linearization and monotonicity methods combining with localized potentials.

\section{The main results}\label{Sec:2}
In this section, we will formulate our two main results: The semilinear elliptic equation \eqref{Main equation} is uniquely solvable (for sufficiently small Dirichlet data), and the associated DN operator uniquely determines the coefficients
in equation \eqref{Main equation}.
	
Our results will be valid under the following assumptions on the domain and the coefficients.

\begin{assum}\label{assum:pcw_anal}
We assume that $\Omega\subset \R^n$, $n\geq 2$ is a bounded domain with $C^\infty$-smooth boundary $\p \Omega$, and that $\sigma\in L^\infty_+(\Omega)$ is a piecewise analytic function in the sense of \cite[Section~3]{KV1985determining}. 
The function $a:\ \Omega\times \R\to \R$ is assumed to fulfill
\[
a(x,y)=\sum_{k=0}^\infty a_k(x) \frac{y^k}{k!},
\]
with
$a_k\in L^\infty(\Omega)$, $a_0=a_1=0$, and $\displaystyle\sup_{k\geq 2}\norm{a_k}_{L^\infty(\Omega)}<\infty$. Moreover each function $a_k$ is assumed to be piecewise analytic in the sense of \cite[Sect.~3]{KV1985determining}.
\end{assum}

Note that \cite[Sect.~3]{KV1985determining} implies that two piecewise analytic functions are piecewise analytic with respect to the same partition, and this naturally extends to every finite number of piecewise analytic functions. However, Assumption \ref{assum:pcw_anal} contains infinitely many such functions $\sigma$, $a_k$ ($k\in \N$), and we
do not assume that they are piecewise analytic with respect to the same partition. 

For our solvability result for the forward problem, we use the following solution space
\[
V:=\left\{ v\in H^1(\Omega):\ -\nabla\cdot (\sigma\nabla v)\in L^\infty(\Omega),\ v|_{\partial \Omega}\in C^{\alpha}(\partial \Omega)\right\}.
\]
equipped with the norm
\[
\norm{v}_V:=\norm{v}_{H^1(\Omega)}+\norm{\nabla\cdot (\sigma\nabla v)}_{L^\infty(\Omega)}+\norm{v|_{\partial \Omega}}_{C^{\alpha}(\partial \Omega)}.
\]
Clearly, $V$ is a Banach space. Moreover, by a result of Li and Vogelius \cite[Corollary 7.3]{LV2000gradient}, 
\begin{equation}
V\subseteq H^1(\Omega)\cap L^\infty(\Omega)  \text{ is continuously embedded.}
\end{equation}

\begin{thm}[Local well-posedness of the forward problem]\label{Thm: Well-posedness}
Let $\Omega \subset \R^n$, $\sigma:\ \Omega\to \R$, and $a:\ \Omega\times \R\to \R$ fulfill assumption \ref{assum:pcw_anal}.		
Then there exists $\veps>0$ so that for all
\begin{align}\label{N_epsilon}
	f\in N_{\veps}:=\left\{ \phi\in C^{\alpha}(\p \Omega): \ \norm{\phi}_{C^{\alpha}(\p \Omega)}<\veps \right\},
\end{align}
there exists a solution $u\in V\subseteq H^1(\Omega)\cap L^\infty(\Omega)$ to the Dirichlet problem  
		\begin{align}\label{Dirichlet problem in well-posedness}
			\begin{cases}
			-\nabla \cdot (\sigma \nabla u)+ a(x,u)=0 & \text{ in }\Omega, \\
			u=f & \text{ on }\p \Omega.
			\end{cases}
		\end{align}
		
Moreover, there exists $\delta>0$, so that, for all $f\in N_{\veps}$, the solution is unique in the set 
of all 
\[
H_\delta:=\{ v\in H^1(\Omega)\cap L^\infty(\Omega):\ \norm{v}_{H^1(\Omega)} + \norm{v}_{L^\infty(\Omega)}\leq \delta\},
\]
and that the solution operator
\[
\mathcal S:\ N_\veps\to V\subseteq H^1(\Omega)\cap L^\infty(\Omega),\quad
f\mapsto u, \quad \text{ where $u\in H_\delta$ solves \eqref{Dirichlet problem in well-posedness}} 
\]
is infinitely differentiable.
\end{thm}

Clearly, $V\subseteq L^\infty(\Omega)$ implies that $x\mapsto a(x,u)$ is a $L^\infty(\Omega)$-function for all $u\in V$. Hence, a solution $u\in V$ of \eqref{Dirichlet problem in well-posedness} has well-defined Neumann boundary values.so that we can define the (non-linear) Dirichlet-Neumann-Operator
\[
\Lambda_{\sigma,a}:\ N_\varepsilon\to H^{-1/2}(\p \Omega),
\quad  \Lambda_{\sigma,a}(f):=\left.\sigma\partial_\nu u \right|_{\p\Omega},
\]
where $u$ is the (sufficiently small) solution of \eqref{Dirichlet problem in well-posedness}. 

Let $D \Subset \Omega$ be an open set with $C^\infty$ boundary $\p D$ such that $\Omega \setminus \overline{D}$ is connected. Then the above result also implies local well-posedness of the Dirichlet problem
\begin{align}\label{BVP with obstacle}
	\begin{cases}
		-\nabla \cdot (\sigma \nabla u) +a(x,u)=0 & \text{ in }\Omega \setminus \overline{D}, \\
		u=0 & \text{ on }\p D,\\
		u=f & \text{ on }\p \Omega,
	\end{cases}
\end{align}
We denote the corresponding DN operator by
\begin{align*}
	\Lambda_{\sigma,a,D}: N_\epsilon \to H^{-1/2}(\p\Omega), \qquad  \Lambda_{\sigma,a,D}	(f) := \left. \sigma \p_\nu u_f\right|_{\p\Omega},
\end{align*}
where $u$ is the (sufficiently small) solution \eqref{BVP with obstacle}.
For an open boundary part $\Gamma\subseteq \p\Omega$, the restriction of $\Lambda_{\sigma,a}$, resp., $\Lambda_{\sigma,a,D}$, to $\Gamma$ is denoted by $\Lambda^\Gamma_{\sigma,a}$, resp., $\Lambda^\Gamma_{\sigma,a,D}$.

Note that the well-posedness of \eqref{BVP with obstacle} can be found in \cite[Appendix]{KU2019derivative_partial} and \cite[Proposition 2.1]{LLLS2019nonlinear} in a slightly different settings that required the coefficients to be sufficiently smooth. However, in this paper, we assume piecewise analytic coefficients so that the coefficients may have jumps. The following theorem extends the uniqueness results from \cite[Theorem 1.2]{LLLS2019partial} and \cite[Theorem 1.6]{KU2019derivative_partial} to this setting.

\begin{thm}[Simultaneous recovering of coefficients and obstacle]\label{Thm: Simultaneous recovery}
	Let $\Omega \subset \R^n$, and two set of coefficients $(\sigma,a)$, and $(\tilde \sigma, \tilde a)$ each  fulfill Assumption \ref{assum:pcw_anal} in connected sets $\Omega \setminus \overline{D}$ and $\Omega \setminus \overline{\tilde D}$, respectively, where $D,\tilde D\Subset \Omega$ are open (possibly empty) sets. Let $\Gamma\subseteq \p\Omega$ be an open boundary part, and let $\varepsilon>0$ be sufficiently small, so that both, $\Lambda^{\Gamma}_{\sigma,a,D}$ and $\Lambda^\Gamma_{\tilde \sigma,\tilde a, \tilde D}$, are defined on $N_\varepsilon$.
	Suppose that 
		\[
		\Lambda^\Gamma_{\sigma,a,D}(f)=\Lambda^\Gamma_{\tilde \sigma,\tilde a,\tilde D}(f)
		\]
		for all $f\in N_\varepsilon$ with $\mathrm{supp}(f)\subseteq \Gamma$, then
		\[
		\sigma=\tilde \sigma, \quad a=\tilde a  \quad\text{ and } \quad  D=\tilde D.
		\]
\end{thm}
To our best knowledge, the preceding theorem is a new result.
The proof will be given in the next section. 

\begin{rmk}Let us emphasize that:
	\begin{itemize}
		\item[(a)] 	Note that even for the full data case, that is, when $\Gamma =\p \Omega$, Theorem \ref{Thm: Simultaneous recovery} is also a new result to our best knowledge. Furthermore, for the case $\Gamma=\p \Omega$, we can reduce the regularity of $\sigma$, by using the first linearization and results in \cite{haberman2013uniqueness} when $\sigma$ is Lipschitz continuous for $n\geq 3$, and \cite{AP2006calderon} when $\sigma\in L^\infty(\Omega)$ for $n=2$. 
		
		\item[(b)] Note that Theorem \ref{Thm: Simultaneous recovery} also covers the case where one of the sets $D$ or $\tilde D$ is empty. Hence, the DN operator also uniquely determines whether there is a cavity or not.
	\end{itemize}
\end{rmk}

\section{Proofs of main results}\label{Sec:3}

To prove our two main results, we start with the following lemma, which will be used for our results.

\begin{lem}\label{lem:well_posedness_V}
The mapping
\[
G:\ V\to L^\infty(\Omega), \quad v(x) \mapsto a(x,v(x))
\]
is infinitely differentiable and its $l$-th Frech\'et derivative fulfills 
\begin{equation*}
G^{(l)}(v)(w_1,\ldots,w_l)=\sum_{k=0}^\infty a_{k+l}(x) \frac{v(x)^{k}}{k!} w_1\ldots w_l \quad \text{ for all } v,w\in V.
\end{equation*}
\end{lem}
\begin{proof}
We also define for $l\in \N_0$
\[
G_l:\ V\to L^\infty(\Omega), \quad v(x)\mapsto \sum_{k=0}^\infty a_{k+l}(x)\frac{v(x)^{k}}{k!}
\]
Then $G_0=G$, and for all $v\in V$, $l\in \N_0$, $G_l(v)\in L^\infty(\Omega)$ follows from $V\subseteq L^\infty(\Omega)$. 

We will prove that $G_l$ is one-time Frech\'et differentiable and that its derivative $G_l':\ V\to \mathcal L(V,L^\infty(\Omega))$ is given by 
\begin{equation}\label{eq:well_posedness_V_Gprime}
G_l'(v)w=M_w G_{l+1}(v) \quad \text{ for all } v,w\in V,
\end{equation}
where, for $w\in V\subseteq L^\infty(\Omega)$
\[
M_w:\ L^\infty(\Omega)\to L^\infty(\Omega),\quad u\mapsto w u
\]
denotes the continuous linear multiplication operator. Then the assertion follows by trivial induction.

Clearly, \eqref{eq:well_posedness_V_Gprime} defines a continuous linear operator $G_l'(v)\in \mathcal L(V,L^\infty(\Omega))$.
To prove that this is indeed the Fr\'echet derivative of $G_l$, let
$v,w\in V$, $x\in \Omega$, and define 
\begin{align*}
	\psi_x&:\ \R\to \R, \\
	\psi_x(t)&:=G_l(v+t(w-v))(x)=\sum_{k=0}^\infty a_{k+l}(x) \frac{(v(x)+t(w(x)-v(x))^k}{k!}.
\end{align*}
Then $\psi$ is infinitely differentiable with 
\begin{align*}
\psi_x'(t)&=\sum_{k=1}^\infty a_{k+l}(x) \frac{(v(x)+t(w(x)-v(x))^{k-1}}{(k-1)!}(w(x)-v(x))\\
\psi_x''(t)&=\sum_{k=2}^\infty a_{k+l}(x) \frac{(v(x)+t(w(x)-v(x))^{k-2}}{(k-2)!}(w(x)-v(x))^2.
\end{align*}
Using 
\[
|v(x)+t(w(x)-v(x))|\leq |v(x)|+|w(x)| \quad \text{ for all } t\in [0,1],
\]
and that $a_k$ are uniformly bounded, and that $V\subseteq L^\infty(\Omega)$ is continuously embedded 
we have that 
\[
\left|\psi_x''(t) \right|\leq C \norm{w-v}_{V}^2 \exp\left(\norm{v}_V+\norm{w}_V\right) \quad \text{ for all } x\in \Omega.
\]
Using Taylor's formula
\[
\left|\psi_x(1)-\psi_x(0)-\psi_x'(0) \right|\leq \frac{1}{2}\max_{\tau\in [0,1]} |\psi''(\tau)|,
\]
we thus obtain 
\[
\norm{G_l(v)-G_l(w)-G_l'(v)(w-v)}_{L^\infty(\Omega)}
\leq \frac{C}{2} \norm{w-v}_{V}^2 \exp\left(\norm{v}_V+\norm{w}_V\right),
\]
so that the assertion is proven.
\end{proof}

\subsection{Local well-posedness result for the forward problem.}

\begin{proof}[Proof of Theorem \ref{Thm: Well-posedness}]
We will apply the implicit function theorem to the map 
\begin{align*}
F:\ C^{\alpha}(\partial \Omega)\times V & \to W:=L^\infty(\Omega) \times C^{\alpha}(\partial \Omega),\\
F:\ (f,v)&\mapsto  \left( -\tre{\nabla\cdot}(\sigma\nabla v)+a(x,v), v|_{\partial \Omega}-f\right).
\end{align*}
By Lemma~\ref{lem:well_posedness_V} this mapping is well-defined and infinitely differentiable, 
and its derivative with respect to $v\in V$ is the continuous linear operator
\[
D_v F(0,0):\ V\to W 
\]
with
\[
D_v F(0,0) u= \left( -(\sigma\nabla u), u|_{\partial \Omega}\right).
\]
Given $w=(w_1,w_2)\in W$ there exists a unique solution $u\in H^1$ of
\[
-\tre{\nabla \cdot}(\sigma\nabla u)=w_1,\quad \text{ and } \quad u|_{\partial \Omega}=w_2.
\]
Then $u\in V$ holds by definition of $V$, which shows that $D_v F(0,0)$ is surjective. 
Since $u$ is the unique solution, $D_v F(0,0)$ is also injective. Since also $F(0,0)=0$ is fulfilled, we can apply the
implicit function theorem (cf., e.g., \cite[Sect.~10.1.1]{renardy2006introduction}), which yields 
an infinitely differentiable function
\[
S:\ N_{\veps}\to V
\]
defined on a neighborhood of the origin $C^{\alpha}(\p \Omega)$,
\[
N_{\veps}:=\{ \phi\in C^{\alpha}(\p \Omega): \ \norm{\phi}_{C^{\alpha}(\p \Omega)}<\veps \},
\]
so that
\[
F(f,S(f))=0 \quad \text{ for all } f\in N_\veps,
\]
and $S(f)$ is the only such element in a neighborhood of the origin in $V$.
Since $F(f,S(f))=0$ implies that $S(f)\in V$ solves \eqref{Dirichlet problem in well-posedness},
the existence of a solution in $V\subseteq H^1(\Omega)\cap L^\infty(\Omega)$ is proven.
Moreover, since every solution $u\in H^1(\Omega)\cap L^\infty(\Omega)$ of \eqref{Dirichlet problem in well-posedness}
fulfills
\[
u\in V,\quad \text{ and } \quad \norm{u}_V\leq \norm{u}_{H^1(\Omega)}
+ \sup_{k\geq 2}\norm{a_k}_{L^\infty(\Omega)} e^{\norm{u}_{L^\infty(\Omega)}}+
\norm{f}_{C^\alpha(\p \Omega)},
\]
the solution is unique in $H_\delta$ with sufficiently small $\delta$. 
\end{proof}

\subsection{The higher order linearization}

To prove Theorem~\ref{Thm: Simultaneous recovery} we first derive some auxilliary results on the
higher-order derivatives of the solution of \eqref{Dirichlet problem in well-posedness}. In the rest of this note, let us fix $\varepsilon >0$ to be a sufficiently small number, such that the well-posedness for \eqref{Main equation} and \eqref{Main equation_obstacle} hold, for any $f\in N_\varepsilon$.

\begin{lem}\label{lemma:derivatives_solution}
Let $f_1, f_2\in N_\varepsilon$, and define
\[
F:\  \LC -1/2, 1/2\RC \times \LC -1/2, 1/2\RC \to V, \quad F(t_1,t_2):=\mathcal S(t_1f_1 + t_2 f_2).
\]
Then
\[
u_1^{(\ell)} :=\p_{t_\ell}F(0,0)=\left. \partial_{t_\ell} \mathcal S(t_1f_1 +t_2f_2)\right|_{t_1=t_2=0}\in V
\] 
solves
\begin{align*}
\begin{cases}
	\nabla \cdot (\sigma \nabla u_1^{(\ell)})=0 &\text{ in $\Omega$}\\
	u_1^{(\ell)}=f_\ell & \text{ on $\partial \Omega$},
\end{cases}
\end{align*}
for $\ell=1,2$. Moreover, for all $m>1$, $m\in \N$, 
\[
u_m:=\p_{t_1}^2 \p_{t_2}^{m-2}F(0,0)=\left. \partial_{t_1}^2 \p_{t_2}^{m-2} \mathcal S(t_1f_1+t_2f_2)\right|_{t_1=t_2=0}\in V
\] 
solves 
\begin{align*}
\begin{cases}
	\nabla \cdot (\sigma \nabla u_m)=\left. \partial_{t_1}^2 \p_{t_2}^{m-2} G(F(t))\right|_{t_1=t_2=0} & \text{ in $\Omega$}\\
	u_m=0  & \text{ on $\partial \Omega$.}
\end{cases}
\end{align*}
Moreover, for all $f\in N_\varepsilon$, the mapping
\[
(t_1,t_2)\mapsto \Lambda(t_1f_1 + t_2f_2),\quad \LC -1/2, 1/2\RC \times \LC -1/2,1/2 \RC \to H^{-1/2}(\p\Omega)
\]
is infinitely differentiable, and
\[
\partial_{t_1}^2 \p _{t_2}^{m-2} \Lambda(t_1f_1 + t_2 f_2)|_{t_1 =t_2=0}=\left. \sigma\partial_\nu u_m\right|_{\p\Omega}.
\]
\end{lem}
\begin{proof}
Note that $F$, $G$, and $\mathcal S$ are infinitely differentiable functions by Lemma~\ref{lem:well_posedness_V}, and Theorem \ref{Thm: Well-posedness}. Moreover, since the trace operator $u\mapsto u|_{\partial \Omega}$ is a continuous linear function from $V$ to $H^{1/2}(\Omega)$,
it follows that
\[
\left. u_1^{(\ell)}\right|_{\partial \Omega}=f_\ell, \text{ for }\ell=1,2, \quad \text{ and } \quad  \left. u_m \right|_{\partial \Omega}=0,  \text{ for $ m\geq 2$.}
\]
Let $v\in H_0^1(\Omega)$. Since $u:=\mathcal S(t_1f_1 + t_2f_2)$ solves \eqref{Dirichlet problem in well-posedness} we have
that
\begin{align*}
0&=\int_{\Omega} \sigma \nabla u \cdot \nabla v \dx + \int_\Omega a(x,u)v  \dx\\
&=\int_{\Omega} \sigma \nabla \mathcal S(t_1f_1+t_2f_2) \cdot \nabla v \dx + \int_\Omega G(\mathcal S(t_1f_1+t_2f_2))v  \dx.
\end{align*}
Noting that the first summand is a linear continuous functional with respect to $\mathcal S(t_1f_1 + t_2f_2)\in V$, 
and the second summand is linear and continuous with respect to $G(\mathcal S(t_1f_1 + t_2f_2))\in L^\infty(\Omega)$, we obtain by differentiation
\begin{align*}
0=&\int_{\Omega} \left. \sigma \nabla \partial_{t_1}^2 \p_{t_2}^{m-2} \mathcal S(t_1f_1 + t_2f_2)\right|_{t_1=t_2=0} \cdot \nabla v \dx \\
&+ \int_\Omega \left. \partial_{t_1}^2 \p_{t_2}^{m-2} G(\mathcal S(t_1f_1 + t_2f_2))\right|_{t_1=t_2=0} v \dx\\
=&\int_{\Omega} \sigma \nabla u_m \cdot \nabla v \dx + \int_\Omega \left. \partial_{t_1}^2 \p_{t_2}^{m-2} G(\mathcal S(t_1f_1 + t_2f_2)\right|_{t_1=t_2=0} v \dx.
\end{align*}
This proves that, for all $m > 1$,
\[
\nabla \cdot (\sigma \nabla u_m)=\left. \partial_{t_1}^2 \p_{t_2}^{m-2} G(\mathcal S(t_1f_1 + t_2f_2))\right|_{t_1=t_2=0}.
\]
On the other hand, for $\ell=1,2$, it follows from $S(0)=0$, and $G'(0)=0$, that
\[
\left. \partial_{t_\ell}G(\mathcal S(t_1f_1 +t_2 f_2))\right|_{t_1=t_2=0}=G'(\mathcal S(0))\LC \p_{t_\ell}\mathcal S(0)\RC f_\ell=0.
\]
Hence, $u_1^{(\ell)}\in V$ solves $\nabla \cdot (\sigma \nabla u_l^{(\ell)})=0$, for $\ell=1,2$.

Moreover, for all $f_1, f_2\in N_\varepsilon$, and $(t_1,t_2)\in (-1/2,1/2)\times (-1/2, 1/2)$, the Neumann data $\Lambda_{\sigma,a}(t_1f_1+t_2f_2)\in H^{-1/2}(\p\Omega)$ fulfills
\[
\langle \Lambda_{\sigma,a}(t_1f_1+t_2f_2),g\rangle=\int_\Omega \sigma \nabla F(t_1,t_2)\cdot \nabla v\dx
+ \int_\Omega G(F(t_1,t_2)) v(x)\dx
\]
for all $v\in H^1(\Omega)$. As above, the first summand is a linear continuous functional with respect to $F(t_1,t_2)\in V$, and the second summand is linear and continuous with respect to $G(F(t_1,t_2))\in L^\infty(\Omega)$. Hence, $(t_1,t_2)\mapsto \Lambda(t_1f_1+t_2f_2)$ is infinitely differentiable with respect to $t_1$ and $t_2$, and, for all $g\in H^{1/2}(\p\Omega)$,
\begin{align*}
&\left\langle \left.\partial_{t_1}^2 \p_{t_2}^{m-2} \Lambda_{\sigma,a}(t_1f_1 + t_2 f_2)\right|_{t_1 =t_2=0},g \right\rangle\\
=&\left. \partial_{t_1}^2 \p_{t_2}^{m-2} \left\langle \Lambda_{\sigma,a}(t_1f_1 +t_2 f_2),g \right\rangle
\right|_{t_1=t_2=0}\\
=&\int_\Omega \sigma \nabla u_m \cdot \nabla v\dx + \int_\Omega \left. \partial_{t_1}^2 \p_{t_2}^{m-2} G(F(t_1,t_2))\right|_{t_1=t_2=0} v(x)\dx\\
=& \left\langle \sigma \partial_\nu u_m|_{\p\Omega},g\right\rangle,
\end{align*}
which proves that $\left.\partial_{t_1}^2 \p_{t_2}^{m-2} \Lambda_{\sigma,a}(t_1f_1 +t_2 f_2)\right|_{t_1=t_2=0}=\left.\partial_\nu u_m\right|_{\p\Omega}$ as desired.
\end{proof}

\begin{lem}\label{lemma:chain_rule}
Let $m\in \N$. There exist numbers $a_{m,j}^{p_1,p_1',\ldots,p_j,p_j'}\in \N_0$ (depending on $j=1,\ldots,m$, and $p_{1},p_1',\ldots p_{j},p_j'\in \N\cup \{0\}$ with $p_1+p_1'+\ldots+p_j+p_j'=m$) so that
\begin{align*}
&\p_{t_1}^2 \p_{t_2}^{m-2}G(F(t_1,t_2))\\
=&G^{(m)}(F(t_1,t_2))\underbrace{\left(\p_{t_1}F(t_1,t_2),\p_{t_1}F(t_1,t_2), \p_{t_2}F(t_1,t_2), \ldots,\p_{t_2}F(t_1,t_2)\right)}_{m\text{-tuples}}\\
& + \sum_{j=2}^{m-1} \sum_{\mycom{p_1,p_1',\ldots,p_j,p_j'\in \N\cup\{0\},}{ p_1+\ldots,+p_j= 2,\ p_1'+p_2'+\ldots+p_j'=m-2}} a_{m,j}^{p_1,p_1',\ldots,p_j,p_j'} G^{(j)}(F(t_1,t_2)) \\
&\qquad \quad  \times \underbrace{ \left( \p_{t_1}^{p_1}\p_{t_2}^{p_1'}F(t_1,t_2),\ldots, \p_{t_1}^{p_j}\p_{t_2}^{p'_{j}}F(t_1,t_2)\right)}_{j\text{-tuples}}\\
&  + G'(F(t_1,t_2))\left(\p_{t_1}^2 \p_{t_2}^{m-2}F(t_1,t_2)\right).
\end{align*}
\end{lem}
\begin{proof}
This follows by induction using the chain rule for the Fr\'echet derivative.
\end{proof}

We also need the following variant of the localized potentials result in \cite{gebauer2008localized}:
\begin{lem}\label{lemma:loc_pot}
Let $D_1,D_2$ be  two disjoint non-empty sets, where
$D_1\subseteq \Omega$ is open, $D_2\subseteq \overline \Omega$ is closed, 
$\Omega\setminus D_2$ is connected, and $\Gamma \cap \overline \Omega\setminus D_2\neq \emptyset$.Then there exists a sequence $(\phi_k)_{k\in \N}\subset C^\alpha(\partial \Omega)$ with $\supp(\phi_k)\subseteq \Gamma$, and 
\[
\int_{D_2} |v_k|^2 \dx \to 0, \quad \text{ and } \quad \int_{D_1} |v_k|^2 \dx \to \infty,
\]
\tre{where $v_k\in H^1(\Omega)$ is the solution of 
\begin{align*}
	\begin{cases}
		\nabla \cdot (\sigma \nabla v_k)=0 &\text{ in }\Omega, \\
		v_k=\phi_k &\text{ on }\p \Omega,
	\end{cases}
\end{align*}
for $k\in \N$.}
\end{lem}
\begin{proof}
Let $C_0^\alpha(\Gamma)$ denote the closure of the space of all $\phi\in C^\alpha(\partial \Omega)$
with $\mathrm{\phi}\subseteq \Gamma$ with respect to the $C^\alpha(\Gamma)$-norm. 
For $j=1,2$, we define $A_j\in \mathcal L(C_0^\alpha(\Gamma),L^2(D_j))$ by
$A_j:\ \phi \mapsto v|_{D_j}$, where $v\in V$ solves
\begin{align*}
 \begin{cases}
 	\nabla\cdot(\sigma\nabla v)=0  & \text{ in }\Omega,\\
 	v=\phi  &\text{ on }\p \Omega.
 \end{cases}
\end{align*}
Note that $V$ is continuously embedded in $L^2(D_j)$, so that $A_j$ are well-defined.

The assertion is proven if we can show that
\[
\not\exists C>0: \ \norm{A_1 \phi}_{L^2(D_1)}\leq C \norm{A_2 \phi}_{L^2(D_2)},
\]
and, by \cite[Lemma~2.5]{gebauer2008localized}, this is equivalent to proving
\[
\mathcal R(A_1')\not\subseteq \mathcal R(A_2').
\]
The operators $A_j'$ are easily checked to map a source term $\psi\in L^2(D_j)$ to 
the Neumann boundary values of the solution of $\nabla\cdot (\sigma\nabla w)=\psi$ with zero Dirichlet data.
By a standard unique continuation argument, it then follows that $0=R(A_1')\cap\mathcal R(A_2')$, which proves the assertion.
\end{proof}

	\subsection{Unique identifiability result for the inverse obstacle problem.}

Now we can prove our simultaneously unique identifiability result for the inverse coefficient problem.

\begin{proof}[Proof of Theorem \ref{Thm: Simultaneous recovery}.]
	\tre{We first show the case as $D=\tilde D=\emptyset$.}	
Let $\Omega \subset \R^n$, and two set of coefficients $(\sigma,a)$, and $(\tilde \sigma, \tilde a)$ each  fulfill Assumption \ref{assum:pcw_anal}. Let $\epsilon>0$ be sufficiently small, so that both, $\Lambda_{\sigma,a}$ and $\Lambda_{\tilde \sigma,\tilde a}$, are defined on $N_\varepsilon$.

Given $f\in N_\varepsilon$, we define the operators $F$, $G$, $\mathcal S$, and the functions 
$u_m\in V$ ($m\in \N$) as in Lemma~\ref{lemma:derivatives_solution}, Lemma~\ref{lem:well_posedness_V}, and Theorem~\ref{Thm: Well-posedness} using the coefficient pair $(\sigma,a)$. The corresponding entities with $(\sigma,a)$ replaced by $(\tilde \sigma,\tilde a)$ will be 
denoted by \tre{$\tilde F$, $\tilde G$, $\tilde{\mathcal S}$, and $\tilde u_m$} ($m\in \N$). 

We will show that
\begin{enumerate}
\item[(a)] If 
\begin{equation}\label{eq:main_proof2_induction_start}
\left. \left. \partial_{t_\ell} \Lambda^\Gamma_{\sigma,a}(t_1f_1 +t_2f_2)\right|_\Gamma \right|_{t_1=t_2=0}=\left.\left.\partial_{t_\ell} \Lambda^\Gamma_{\tilde \sigma,\tilde a}(t_1f_1 +t_2 f_2)\right|_{\Gamma}\right|_{t_1=t_2=0},
\end{equation}
for all $f_1 , f_2\in N_\varepsilon$, then $\sigma=\tilde \sigma$, and $u_1^{(\ell)}=\tilde u_1^{(\ell)}$, for $\ell=1,2$.
\item[(b)] If, for some $m\geq 2$, $\sigma=\tilde \sigma$, $a_j=\tilde a_j$, $u_j=\tilde u_j$ for all $j=1,\ldots,m-1$, and
\[
\left.\left.\partial_{t_1}^2 \p _{t_2}^{m-2} \Lambda^\Gamma_{\sigma,a}(t_1f_1 + t_2 f_2)\right|_\Gamma\right|_{t_1 =t_2=0}=\left. \left. \partial_{t_1}^2 \p_{t_2}^{m-2} \Lambda^\Gamma_{\tilde \sigma,\tilde a}(t_1f_1 +t_2 f_2)\right|_{\Gamma}\right|_{t_1=t_2=0},
\]
for all $f_1, f_2\in N_\varepsilon$,
then $a_m=\tilde a_m$, and $u_m=\tilde u_m$.
\end{enumerate}
Clearly, this proves Theorem \ref{Thm: Simultaneous recovery} by induction, since we assumed that $a_0=0=\tilde a_0$ and $a_1=0=\tilde a_1$.

To show (a), note that \eqref{eq:main_proof2_induction_start} implies that the local DN operator for the linear elliptic equation $\nabla\cdot (\sigma \nabla v)=0$ is the same 
for the two coefficients $\sigma$ and $\tilde \sigma$. This implies $\sigma=\tilde \sigma$ by the classical Kohn-Vogelius result \cite{KV1985determining}, and the uniqueness of solutions yield that $u_1^{(\ell)}=\tilde u_1^{(\ell)}$ in $\Omega$, for $\ell=1,2$.

To prove (b), note that 
\[
G'(F(0,0)) (\p_{t_1}^2 \p_{t_2}^{m-2}F(0,0))=0=\tilde G'(\tilde F(0,0))(\p_{t_1}^2 \p_{t_2}^{m-2}\tilde F(0,0))
\]
since $ F(0,0)=0=\tilde F(0,0)$ and $ G'(0)=0=\tilde G'(0)$. 
Moreover, 
\[
\p_{t_1}^{p_1}\p_{t_2}^{p_2}F(0,0)=u_p=\tilde u_p=\p_{t_1}^{p_1}\p_{t_2}^{p_2}\tilde F(0,0), \quad \text{ for all }\quad p_1+p_2=1,\ldots,m-1,
\] 
and, for all $w_1,\ldots,w_j$,
\begin{align*}
G^{(j)}(F(0,0))(w_1,\ldots,w_l)=&a_j(x)w_1\ldots w_l\\
=&\tilde a_j(x)w_1\ldots w_l =\tilde G^{(j)}(\tilde F(0))(w_1,\ldots,w_l).
\end{align*}
Using Lemma~\ref{lemma:chain_rule} we thus obtain
\begin{align*}
 \lefteqn{\left.\p_{t_1}^2 \p_{t_2}^{m-2}G(F(t_1,t_2))\right|_{t_1=t_2=0}-\left.\p_{t_1}^2 \p_{t_2}^{m-2}\tilde G(\tilde F(t_1,t_2))\right|_{t_1=t_2=0}}\\
=& G^{(m)}(0)\left. \left(\p_{t_1}F(t_1,t_2),\p_{t_1}F(t_1,t_2), \p_{t_2}F(t_1,t_2), \ldots,\p_{t_2}F(t_1,t_2)\right)\right|_{t_1=t_2=0} \\
&-\tilde G^{(m)}(0)\left. \left(\p_{t_1}\tilde F(t_1,t_2),\p_{t_1}\tilde F(t_1,t_2), \p_{t_2}\tilde F(t_1,t_2), \ldots,\p_{t_2}\tilde F(t_1,t_2)\right)\right|_{t_1=t_2=0}\\
=&a_m u_1^m - \tilde a_m \tilde u_1^m = (a_m-\tilde a_m)u_1^m.
\end{align*}
Hence, $u_m-\tilde u_m\in V$ solves
\begin{align}
 \begin{split}
 	0=&-\nabla \cdot (\sigma\nabla u_m)+\left. \p_{t_1}^2 \p_{t_2}^{m-2} G(F(t))\right|_{
 		t_1=t_2=0} \\
 	 &+ \nabla \cdot (\sigma\nabla \tilde u_m) - \left.\p_{t_1}^2 \p_{t_2}^{m-2} G(\tilde F(t))\right|_{t_1=t_2=0}\\
 	=& -\nabla \cdot (\sigma\nabla (u_m-\tilde u_m))+ (a_m-\tilde a_m) \LC  u_1^{(1)}\RC^2 \LC u_1^{(2)} \RC^{m-2}.\label{eq:um_minus_tildeum}
 \end{split}
\end{align}
For all $g\in C^\alpha(\Gamma)$ with $\supp(g)\subset \Gamma$, we thus obtain
\begin{align}	\label{eq:am_minus_tileam_var}
\begin{split}
	0=&\left\langle \left. \left. \p_{t_1}^2 \p_{t_2}^{m-2} \Lambda_{\sigma,a}(t_1f_1+t_2f_2)\right|_\Gamma\right|_{t_1=t_2=0},g\right\rangle \\
	&-\left\langle \left. \left. \p_{t_1}^2 \p_{t_2}^{m-2} \Lambda_{\tilde \sigma,\tilde a}(t_1f_1+t_2f_2)\right|_\Gamma\right|_{t_1=t_2=0},g\right\rangle \\
	 =&\langle \sigma \partial_\nu (u_m-\tilde u_m)|_{\Gamma},g\rangle- \langle \sigma \partial_\nu \tilde u_m|_{\Gamma},g\rangle\\
	 =&\int_\Omega \sigma \nabla (u_m-\tilde u_m)\cdot \nabla v_1\dx + \int_\Omega (a_m-\tilde a_m)\LC u_1^{(1)} \RC^2 \LC  u_1^{(2)}\RC^{m-2}v_1 \dx\\
	 =&\int_\Omega (a_m-\tilde a_m)\LC u_1^{(1)} \RC^2 \LC  u_1^{(2)}\RC^{m-2} v_1\dx,
\end{split}
\end{align}
where $v_1\in V$ solves
\begin{align*}
\begin{cases}
	\nabla \cdot (\sigma \nabla v_1)=0 &\text{ in $\Omega$,}\\
	v_1=g & \text{ on $\partial \Omega$.}
\end{cases}
\end{align*}
We will now show that this implies $a_m=\tilde a_m$. Clearly this also implies $u_m=\tilde u_m$ by using \eqref{eq:um_minus_tildeum} and $\left.u_m \right|_{\p\Omega}=f=\left.\tilde u_m \right|_{\p\Omega}$.

Assume that this is not the case. Since $a_m$ and $\tilde a_m$ are piecewise analytic, 
we can choose two disjoint non-empty sets $D_1,D_2$, where
$D_1\subseteq \Omega$ is open, $D_2\subseteq \overline \Omega$ is closed, 
$\Omega\setminus D_2$ is connected, and $\Gamma \cap \overline \Omega\setminus D_2\neq \emptyset$ 
such that either 
\begin{enumerate}
\item[(i)] $\left. a_m\right|_{\Omega\setminus D_2}\geq  \left. \tilde a_m\right|_{\Omega\setminus D_2}$ and $\left. (a_m-\tilde a_m)\right|_{D_1}\in L_+^\infty(D_1)$, or
\item[(ii)] $\left. a_m\right|_{\Omega\setminus D_2}\leq \left. \tilde a_m\right|_{\Omega\setminus D_2}$ and $\left. (\tilde a_m- a_m)\right|_{D_1}\in L_+^\infty(D_1)$,
\end{enumerate}
cf.\ \cite[Appendix~A]{harrach2013monotonicity} for a proof for $\Gamma=\p\Omega$ that also holds for arbitrarily small open boundary pieces $\Gamma\subset \p\Omega$. Without loss of generality we will asumme that (i) holds true in the following.

For $m\in \N$, let us choose a non-negative, but not identically zero, function $\psi\in N_\varepsilon\subset C^{\alpha}(\p \Omega)$ with $\mathrm{supp}(\psi)\subseteq \Gamma$. By the strong maximum principle, the corresponding solution $w\in V$ of $\nabla \cdot (\sigma\nabla w)=0$ in $\Omega$, with $w|_{\partial \Omega}=\psi$, will be positive inside $\Omega$. 
Then we can use the localized potentials result in Lemma~\ref{lemma:loc_pot} to obtain a sequence $(\phi_k)_{k\in \N}\subset C^\alpha(\partial \Omega)$ with $\supp(\phi_k)\subseteq \Gamma$, and 
\[
\int_{D_2} w_{1,k}^{2} \dx \to 0, \quad \text{ and } \quad \int_{D_1} w_{1,k}^{2} \dx \to \infty,
\]
where $w_{1,k}\in V$ solves $\nabla \cdot (\sigma\nabla w_{1,k})=0$ in $\Omega$ with $\left. w_{1,k}\right|_{\partial \Omega}=\phi_k$.

Using the Dirichlet data
$$
f_{1,k}:=\frac{\varepsilon}{2\norm{\phi_k}_{C^\alpha(\partial \Omega)}}\phi_k\in N_\varepsilon,\quad  f_{2} =\psi  \quad \text{ and }\quad g_k:=\left(\frac{2\norm{\phi_k}_{C^\alpha(\partial \Omega)}}{\varepsilon}\right)^2 \psi,
$$ 
such that solutions of the first linearized equation satisfy $\left. u_{1,k}\right|_{\p \Omega}=f_{1,k}$, $\left. u_{2}\right|_{\p \Omega}=f_{2}$ and $\left. v_{1,k}\right|_{\p \Omega}=g_k$,
then we can obtain from
\eqref{eq:am_minus_tileam_var} that
\begin{align*}
	0=&\int_\Omega (a_m-\tilde a_m)u_{1,k}^{2} u_{2}^{m-2}v_{1,k}\dx \\
	\geq &\int_{D_2} (a_m-\tilde a_m)u_{1,k}^{2} u_{2}^{m-2}v_{1,k} \dx 
	+ \int_{D_1} (a_m-\tilde a_m)u_{1,k}^{2} u_{2}^{m-2}v_{1,k} \dx\to \infty, 
\end{align*}
as $k\to \infty$.
Here we have used the nonnegative of $\psi$ such that $u_2$ and $v_{1,k}$ are positive for any $k\in \N$.
This contradiction shows that (b) holds.

\medskip

On the other hand, if one of $D$ or $\tilde D$ is a nonempty set, similar to the arguments of the previous case, one can determine that $\sigma= \tilde \sigma$ in $\Omega \setminus (\overline{D\cup \tilde D})$ by applying the boundary determination to piecewise analytic functions.		
Let us denote that $u_1$ and $\tilde u_1$ to be the solution of $\nabla \cdot (\sigma \nabla u_1)=0$ and $\nabla \cdot (\tilde \sigma \nabla \tilde u_1)=0$ in $\Omega$, respectively.
By using the strategy introduced in \cite{LLLS2019partial},	let $G$ be the connected component of \tre{$\Omega\setminus (\overline{D \cup \tilde D})$} whose boundary contains $\p \Omega$. Let $U:=u_1-\tilde u_1$, then $U$ is a solution of 
\begin{align*}
	\begin{cases}
		\nabla \cdot \left(\sigma \nabla U \right) =0 & \text{ in }G, \\
		U=\sigma\p_\nu U=0 & \text{ on }\Gamma,
	\end{cases}
\end{align*}
where we have utilized $\Lambda_{\sigma,a,D}^{\Gamma}(f)=\Lambda_{\tilde \sigma,\tilde a,\tilde D}^{\Gamma}(f)$, for any $f\in N_\varepsilon$. By the unique continuation for  second order elliptic equations, one has that $U \equiv 0$ in $G$. Therefore,
\begin{align}\label{uniqueness of first linearization of cavity}
	u_1 = \tilde u_1 \text{ in }G.
\end{align}

We next  prove $D=\tilde D$ via a contradiction argument. Suppose not, i.e., $D \neq \tilde D\Subset  \Omega$, and assume that $\tilde D\neq \emptyset$. By using \cite[Lemma A.1]{LLLS2019partial} or \cite[Section 4]{KU2019derivative_partial}, without loss of generality, we may assume that there exists a point $x_1$ such that 
\begin{align}\label{connected component condition}
	x_1 \in \p G \cap (\Omega \setminus \overline{D})\cap \p \tilde D.
\end{align}
$\tilde u_1(x_1)=0$ since $x_1 \in \p \tilde D$. By \eqref{uniqueness of first linearization of cavity}, we have $u_1(x_1)=0$. Note that $x_1$ is an interior point of the open set $\Omega\setminus \overline{D}$. Consider the boundary values $u_1 |_{\Gamma}\geq 0$ such that $u_1|_{\Gamma}\not \equiv 0$.
Now, since $u_1(x_1)=0$, by the maximum principle, we have that $u_1\equiv 0$ in the connected open set $\Omega \setminus \overline{D}$, which contradicts with the nonzero boundary condition on $\Gamma$. Therefore, the conclusion $D=\tilde D$ must hold. Furthermore, we have by \eqref{uniqueness of first linearization of cavity} that 
\begin{equation}\label{cavity lins}
	u_1=\tilde u_1 \text{ in } \Omega \setminus \overline{D},
\end{equation}
as desired. 
Finally, by repeating the same arguments as in \tre{the case (1)}, we can show that $a=\tilde a$ in $\LC\Omega \setminus \overline{D}\RC \times \R$ as we wish. This proves the assertion.
\end{proof}

From the proof of Theorem \ref{Thm: Well-posedness}, one can see that if one of $D$ or $\tilde D$ is an empty set, then $D=\tilde D=\emptyset$ immediately.

 \medskip

\noindent\textbf{Acknowledgement.} 
Y.-H. Lin is partially  supported by the Ministry of Science and Technology Taiwan, under the Columbus Program: MOST-110-2636-M-009-007. \tre{The authors want to thank the anonymous reviewer for the careful reading and useful suggestions.}

\bibliographystyle{alpha}
\bibliography{ref}

\end{document}